\documentclass[12pt]{amsart}

\usepackage{latexsym,amssymb,amsfonts,amscd,amsthm,amstext}
\unitlength=1mm\linethickness{0.8pt}

\def\newpic#1{}
\numberwithin{equation}{section}

\newtheorem{theorem}{Theorem}[section]
\newtheorem{proposition}{Proposition}[section]
\theoremstyle{definition}
\newtheorem{definition}{Definition}[section]
\newtheorem{example}{Example}[section]
\newtheorem{remark}{Remark}[section]
\newtheorem{lemma}{Lemma}[section]

\newcommand{\be}{\begin{equation}}

\newcommand{\bt}{\begin{tabular}{c}}
\newcommand{\et}{\end{tabular}}

\begin{document}

\author{Ewa Graczy\'{n}ska}
\address{Opole University of Technology, Institute of Mathematics
\newline
ul. Luboszycka 3, 45-036 Opole, Poland} \email{egracz@po.opole.pl
\hspace{5mm} http://www.egracz.po.opole.pl/}
\author{\fbox{Dietmar Schweigert}}
\address{Technische Universit\"{a}t Kaiserslautern, Fachbereich Mathematik
\newline
Postfach 3049 \\ 67653 Kaiserslautern, Germany}

\title{The dimension of a variety}
\hspace{3cm}
\begin{abstract}
The results of this paper were  presented during the Workshop AAA71
and CYA21 in B\c{e}dlewo, Poland on February 11, 2006.

{\it Derived varieties} were invented by P. Cohn in \cite{1}.

{\it Derived varieties of a given type} were invented by the authors in \cite{4}.

In the paper we deal with derived variety $V_{\sigma}$ of a given variety,
by a fixed hypersubstitution $\sigma$. We introduce the notion of the
{\it dimension of a variety} as the cardinality $\kappa$ of the set of all
proper derived varieties of $V$ included in $V$.

We examine dimensions of some varieties in the lattice of all varieties of a
given type $\tau$. Dimensions of varieties of lattices and  all subvarieties
of regular bands are determined.



\end{abstract}

\maketitle

\section{Notations}
By $\tau$ we denote a fixed type $\tau$: $I \rightarrow N$, where
$I$ is an index set and $N$ is the set of all natural numbers. We
use the definition of an $n$-ary {\it term} of type $\tau$ from
\cite[p. 6]{1}.

$T(\tau)$ denotes the set of all term symbols of type
$\tau$. For a given variety $V$ of type $\tau$, two terms $p$ and $q$ of type $\tau$ are
called {\it equivalent} (in $V$) if the identity $p \approx q$ holds in $V$.

\begin{definition}
For a given type $\tau$, $F$ denotes the set of all fundamental
operations $F = \{f_{i}: i \in I\}$ of type $\tau$, i.e. $\tau(i)$
is the arity of the operation symbol $f_{i}$, for $i \in I$. Let
$\sigma = (t_{i}: i \in I)$ be a fixed choice of terms of type
$\tau$ with $\tau(t_{i}) = \tau(f_{i})$, for every
$i \in I$.\\
Recall from \cite{4} (cf. \cite[p. 13]{1}), that for a given $\sigma$, the extension
of $\sigma$ to the map $\overline{\sigma}$ from the set $T(\tau)$ to $T(\tau)$, leaving all the
variables unchanged and acting on composed terms as:
\begin{center}
$\overline{\sigma}(f_{i}(p_{0},...,p_{n-1})) = \sigma(f_{i})(\overline{\sigma}(p_{0}),...,\overline\sigma(p_{n-1}))$
\end{center}
is called a {\it hypersubstitution} of type $\tau$.

In the sequel, we shall use $\sigma$ instead of $\overline{\sigma}$ for a hypersubstitution.\\
A hypersubstitution $\sigma$ will be called {\it trivial}, if it is the identity
mapping.

The set of all hypersubstitutions of type $\tau$ will be denoted by $H(\tau)$.

For any algebra ${\bf A} = (A, \Omega) = (A, (f_{i}^{\bf A}: i \in I)) \in V$,
of type $\tau$, the algebra ${\bf A}_{\sigma}= (A, (t^{\bf A}_{i}: i \in I))$ or shortly
${\bf A}_{\sigma} = (A, \Omega_{\sigma})$, for $\Omega_{\sigma} = (t_{i} :i \in I)$
is called a {\it derived algebra } (of a given type $\tau$) of ${\bf A}$,
corresponding to $\sigma$, for any $\sigma \in H(\tau)$ (cf. \cite{4}, \cite{10}).
\end{definition}

\begin{definition}
The variety generated by the class of all derived algebras ${\bf A}_{\sigma}$,
of algebras ${\bf A} \in V$ will be called the {\it derived variety} of $V$
using $\sigma$ and it will be denoted by $V_{\sigma}$, for any fixed
$\sigma \in H(\tau)$.

For a class $K$ of algebras of a given type $\tau$, $D(K)$ denotes
the class of all derived algebras of $K$ for all possible choices of
$\sigma$ of type $\tau$, i.e.:
\begin{center}
$D(K) = \bigcup \{ K_{\sigma}: \sigma \in H(\tau)\}$.
\end{center}
\end{definition}
$D$ is a class operator examined in \cite{4} (cf. \cite{9}, \cite{10}).

Let us note, that $V_{\sigma} = HSP(\sigma(V))$, for a given variety $V$ and $\sigma$,
where $\sigma(V)$ denotes the class of all derived algebras ${\bf A}_{\sigma}$,
for ${\bf A} \in V$.

Recall from \cite{5}:
\begin{definition}
For a given set $\Sigma$ of identities of type $\tau$, $E(\Sigma)$
denotes the set of all consequences of $\Sigma$ by the rules (1) -- (5)
of inferences of G. Birkhoff (cf. \cite{14}, \cite{5}).

$Mod(\Sigma)$ denotes the variety of algebras determined by $\Sigma$.

A variety $V$ is {\it trivial} if all algebras in $V$ are {\it trivial}
(i.e. one-element). Trivial varieties will be denoted by $T$. A subclass $W$
of a variety $V$ which is also a variety is called {\it subvariety} of $V$.

$V$ is a {\it minimal} (or {\it equationally complete})
variety if $V$ is not trivial but the only subvariety of $V$, which is not
equal to $V$ is trivial.
\end{definition}

We accept the following definition from \cite{10}:

\begin{definition}
A derived variety $V_{\sigma}$ is {\it proper} if $V_{\sigma}$ is not
equal to $V$, i.e. $V_{\sigma} \neq V$.
\end{definition}
Note, that $V_{\sigma}$ may be not proper only for nontrivial $\sigma$.

Recall from \cite{4}:
\begin{definition}
A variety $V$ of type $\tau$ is {\it solid} if $V$ contains all derived
varieties $V_{\sigma}$ for every choice of $\sigma$ of type $\tau$,
i.e. $D(V) \subseteq V$.
\end{definition}
\begin{definition}
A variety $V$ of type $\tau$ is {\it fluid} if the variety $V$ contains
no proper derived varieties $V_{\sigma}$ for every choice of $\sigma$ of
type $\tau$.
\end{definition}
Fluid varieties appear naturally in many well known examples (cf. \cite{12}).
Derived varieties are an important tool for describing the lattice of all
subvarieties of a given variety and therefore we expect some practical
applications of the invented notion.

Note, that our definition of a {\it fluid variety} does not coincide with
that of \cite{10}.
\section{The Dimension}
\begin{definition}
If $V$ is a variety of type $\tau$, then the dimension of $V$ is the
cardinality $\kappa$ of the set of all proper derived varieties $V_{\sigma}$
of $V$ included in $V$, for $\sigma \in H(\tau)$. We write then that $\kappa = dim(V)$.
\end{definition}

From the definitions above it follows that the trivial variety $T$ of a given
type is of dimension 0.

\begin{theorem}
Minimal varieties  are of dimension 0. Fluid varieties are of
dimension 0.
\end{theorem}

Later on we shall use the well-known {\it conjugate property} of \cite{13}
(cf. \cite[p. 35]{egracz} and \cite{12}) and quote as:

\begin{theorem}
Let ${\bf A}$ be an algebra and  $\sigma$ be a hypersubstitution of type $\tau$.
Then an identity $p \approx q$ of type $\tau$ is satisfied in the derived
algebra ${\bf A}_{\sigma}$ if and only if the derived identity
$\sigma(p) \approx \sigma(q)$ holds in ${\bf A}$.
\end{theorem}
From the theorem above, it immediately follows:
\begin{theorem}
Let  $V$ be a variety and two hypersubstitutions $\sigma_{1}$ and
$\sigma_{2}$ of type $\tau$ be given.
If $\sigma_{1}(f_{i}) \approx \sigma_{2}(f_{i})$, is an identity of $V$ for
every $i \in I$, then the derived varieties $V_{\sigma_{1}}$ and
$V_{\sigma_{2}}$ are equal.
\end{theorem}
{\it Proof}. The proof follows by induction on the complexity of terms of
type $\tau$. $\Box$

In the proof we use the relation $\sim_{V}$ on sets of
hypersubstitutions which was introduced by J. P{\l}onka in
\cite{JP1} and used in \cite{13} to determine the notion of $V$-{\it
equivalent hypersubstitutions} in order to simplify the procedure of
checking whether an identity is satisfied in a variety $V$ as a
hyperidentity.

\begin{theorem}
Assume that in a variety $V$ (of a finite type) there is only a finite
number of non-equivalent n-ary terms, for every $n \in N$. Then $V$ is
of a finite dimension.
\end{theorem}
{\it Proof}. The proof follows from the fact that under the assumption,
in $V$ there are only finitely many non-equivalent fundamental operations
and  non-equivalent hypersubstitutions of type $\tau$. Therefore there are
only finitely many derived varieties of $V$ and $dim(V)$ is finite. $\Box$

\section{Dimensions of varieties of lattices}

We  present some examples in lattice varieties as an answer to
a problem posed by Brian Davey (La Trobe University, Australia) during
the Conference on Universal Algebra and Lattice Theory (July 2005)
at Szeged University (Hungary).

Let ${\bf L} = (L, \vee, \wedge)$ be a lattice. A variety $L_{\sigma}$
derived from a variety $L$ of lattices must not be a variety of lattices.

This follows from the fact, that there are only four non-equivalent binary
terms in lattices, namely $x$, $y$, $x \vee y$ and $x \wedge y$. Given a
hypersubstitution $\sigma$ of type (2,2). If $\sigma$ is trivial, then the derived algebra
${\bf L}_{\sigma}$ is ${\bf L}$ itself.  If one takes $\sigma$ generated
by $\sigma(\vee) = \wedge$ and $\sigma(\wedge) = \vee$, then ${\bf L}_{\sigma}$ is the dual lattice
${\bf L}^{d} = (L, \wedge, \vee)$. Otherwise the derived algebra
${\bf L}_{\sigma}$ is not a lattice at all, as some lattice axioms will be failed, unless ${\bf L}$ is
trivial (i.e. one-element lattice).

We got immediately:

\begin{example}
Let $V$ be a nontrivial variety of lattices. Then a derived variety
$V_{\sigma}$ is the dual variety of lattices $V^{d}$ or a variety which is
not a variety of lattices.
\end{example}
\begin{example}
The variety $L$ of all lattices in type (2,2) is fluid and not solid.
\end{example}
The variety $L$ is fluid as it is selfdual, i.e. $L = L^{d}$.
It is not solid, as the commutativity laws for $\vee$ and $\wedge$ are not
satisfied as hyperidentities in lattices, for example.
\begin{theorem}
Every variety of lattices is fluid.
\end{theorem}
{\it Proof}. Let $V$ be a variety of lattices. Consider the dual
variety of $V$, i.e. the variety $V^{d}$ of all dual lattices of
$V$. Then there are only two
 possibilities:

(i) $V^{d} \subseteq V$ and consequently $V=V^{d}$

or

(ii) $V$ and $V^{d}$ are incomparable in the lattice of all varieties of lattices.

Therefore we conclude, that either $V$ is selfdual or $V$ and
$V^{d}$ are incomparable. In consequence $V$ is fluid and
$dim(V)=0$.  $\Box$

\section{Dimensions of subvarieties of regular bands}
In this section we concentrate on the lattice of all subvarieties of regular
bands, described in \cite{fen1} -- \cite{fen2} and \cite{3}.

\begin{definition}
{\it Bands} is the variety $B$ of algebras of type (2), defined by:
associativity and idempotency (i.e. a {\bf band} is an idempotent semigroup).
\end{definition}

Following \cite[p. 11]{2},  let us note, that the variety of bands
has only six non-equivalent binary terms, therefore only six
hypersubstitutions of type (2) in the variety of bands should be
checked, namely: $\sigma_{1} - \sigma_{6}$ defined as follows:
$\sigma_{1}(xy)=x$, $\sigma_{2}(xy)=y$, $\sigma_{3}(xy)=xy$,
$\sigma_{4}(x,y) = yx$,  $\sigma_{5}(xy)=xyx$ and
$\sigma_{6}(xy)=yxy$ to be considered in order to determine all
derived varieties of a given subvariety of regular bands.

Recall Proposition 3.1.5(i) from \cite[p. 11, 77]{1}:

\begin{definition}
A variety $V$ of type (2) is called {\it hyperassociative} if the
associativity law is satisfied in $V$ as a hyperidentity.
\end{definition}
\begin{proposition}
A variety of bands is hyperassociative if and only if it is contained in the
variety  $RegB$  of regular bands.
\end{proposition}
The propositions above may be considered as a motivation of our interest in
the lattice of all subvarieties of the variety of regular bands.

In order to determine the dimension of all subvarieties of RegB,
we shall use the following two theorems of \cite{12}:
\begin{theorem}
The variety of $B$ all bands constitutes a not fluid and not solid variety
of type (2).
\end{theorem}

\begin{theorem}
A variety $V$ of bands is fluid if and only if it is minimal.
\end{theorem}
\begin{remark}
Note, that a nontrivial variety $V$ is of dimension 0 if and only if it is fluid.
\end{remark}
\begin{definition}
An identity $e$ of the form $p\approx q$ is called {\it leftmost} ({\it rightmost})
if and only if it has the same first (last) variable on each side.
One which meets both of these conditions is called {\it outermost}.
\end{definition}

First we express three technical lemmas:
\begin{lemma}
Let $\Sigma$ be a set of identities of type $\tau$ which are leftmost
(or rightmost). Then the set $E(\Sigma)$ consists only of leftmost
(rightmost) identities.
\end{lemma}
{\it Proof}. The proof follows from the observation that all rules of
inference (1) -- (5) preserves the property of being the {\it leftmost} (or {\it rightmost})
identity. Therefore the closure of the set of left(right)most identities
consists of left(right)most identities. $\Box$

From \cite{fen1}--\cite{fen2} and \cite{3} it follows that every subvariety of
the variety $B$ of all bands is defined by one additional identity added to
two axioms of bands (i.e. associativity and idempotency).

\begin{lemma}      \label{lemma2}
Assume that $V$ and $W$ are varieties of bands, $W$ is  defined by a single
identity $p\approx q$, i.e. $W = Mod(p\approx q)$ (in the varieties of bands). Then:
\begin{center}
$V_{\sigma} \subseteq W$, for a given $\sigma \in H(\tau)$,
\end{center}
if and only if the derived identity $\sigma(p) \approx \sigma(q)$ is satisfied in
$V$, i.e. $V \models \sigma(p) \approx \sigma(q)$.
\end{lemma}
{\it Proof}. $V_{\sigma} = HSP(\sigma(V)) \subseteq W$ if and only if
$\sigma(V) \models p\approx q$. By theorem 2.2 we conclude that
${\bf A}_{\sigma} \models p\approx q$, for every algebra ${\bf A}_{\sigma} \in \sigma(V)$,
if and only if ${\bf A} \models \sigma(p) \approx \sigma(q)$, for every algebra
${\bf A} \in V$, i.e. $V \models \sigma(p) \approx \sigma(q)$. $\Box$

A simple generalization of the above lemma is the following:

\begin{lemma}
Assume that $V$ and $W$ are varieties of type $\tau$, $W$ is  defined by
a set $\Sigma$ of identities of type $\tau$, i.e. $W = Mod(\Sigma)$. Then:
\begin{center}
$V_{\sigma} \subseteq W$, for a given $\sigma \in H(\tau)$,
\end{center}
if and only if the derived identity $\sigma(p) \approx \sigma(q)$ is satisfied in
$V$, i.e.
\begin{center}
$V \models \sigma(p) \approx \sigma(q)$, for every identity $p \approx  q \in
\Sigma$.
\end{center}
\end{lemma}
{\it Proof}. The proof is similar as that of lemma \ref{lemma2}, where
$p \approx q$ is any identity of the given axiomatic $\Sigma$. $\Box$

The next three propositions show some regularities in the dimensions of
all subvarieties of regular bands described in \cite[p. 244]{fen2} and
\cite{3}:

\begin{definition}
The variety of $RB$ in the variety  $B$ of bands is defined by
the identity: $y \approx  yxy$.  It is called the variety of rectangular bands.
\end{definition}
The fact that the variety $RB$ is solid was proved in \cite[p. 96]{2}.

We expressed the situation of theorems above on the diagram, which
describes the bottom part of the lattice of all identities of bands,
see \cite{4} and \cite[p. 244]{5} Proposition 3.1.5 of \cite{1}:

\vspace{5mm}

\input{dim1.pic}

\begin{theorem}
The variety $RB$ is of dimension 2.
\end{theorem}
{\it Proof}. The variety of $RB$ of rectangular bands has only two
subvarieties, namely the variety $LZ$ defined by the identities:
$yx\approx y$ (called the variety of left-zero semigroups) and the variety
$RZ$ defined by $xy \approx y$ (called the variety of right-zero
semigroups), respectively. Both of them are derived varieties of
$RB$ by the first and the second projection, respectively. To prove
that, let $({\bf A})_{\sigma_{1}} \in \sigma_{1}(RB)$, for ${\bf A}
\in RB$. Then the identity $yx \approx y$ is satisfied in ${\bf
A}_{\sigma_{1}}$, as: $\sigma_{1}(yx) \approx y \approx y \approx \sigma_{1}(y)$ is
satisfied in ${\bf A}$ and consequently in $(RB)_{\sigma_{1}}$.
Similarly for $\sigma_{2}$. We conclude that $dim(RB) = 2$. $\Box$
\begin{theorem}
The varieties $V_{1}$ and $V_{2}$ of bands defined by the identities:
\begin{center}
(1)  $zxy\approx zyx$ \hspace{5mm} and \hspace{5mm} (2) $yxz \approx xyz$, respectively,
\end{center}
are mutually derived by $\sigma_{4}$. Moreover, $dim(V_{1}) = dim(V_{2}) = 1$.
\end{theorem}
{\it Proof}. Note, that the varieties $V_{1}$ and $V_{2}$ has only
two proper subvarieties, namely: the variety of left (right)
zero-semigroups (respectively) and the variety $SL$ of semilattices.
The variety of semilattices, defined (in the variety of bands) by
the commutativity law: $xy\approx yx$ is not a derived variety of
$V_{1}$, nether of $V_{2}$. This follows from the fact, that if the
variety $SL$ of semilattices would be a derived variety of $V_{1}$,
then $SL = (V_{1})_{\sigma_{5}}$ or $SL = (V_{1})_{\sigma_{6}}$.
This is impossible, via theorem 1.3, as the derived identity of
$xy\approx yx$ by the hypersubstitions $\sigma_{5}$ (or
$\sigma_{6}$), i.e. $\sigma_{5}(xy) \approx \sigma_{5}(yx)$ (or
$\sigma_{6}(xy) \approx \sigma_{6}(yx)$) is of the form $xyx\approx
yxy$ is neither leftmost nor rightmost  and therefore, by lemma 4.1
is not satisfied in $V_1$ as every identity satisfied in $V_{1}$ is
leftmost. Similarly for $V_{2}$. The proof follows from the fact
that the only proper derived variety of $V_{1}$ included in $V_{1}$
is the variety $LZ$ of  left zero semigroups. Similarly, one can
show that the only proper derived  subvariety of $V_{2}$ by the
second projection $\sigma_{2}$ is the variety $RZ$ of right zero
semigroups. Finally we conclude that $dim(V_{1})=dim(V_{2})=1$.
$\Box$
\begin{definition}
Varieties of dimension 1 will be called {\it prefluid}.
\end{definition}
\begin{theorem}
The varieties $V_{3}$ and $V_{4}$ of bands defined by the identities:
\begin{center}
(3)  $yx \approx yxy$ \hspace{5mm} and \hspace{5mm} (4) $xy \approx yxy$, respectively,
\end{center}
are mutually derived by $\sigma_{4}$.  Moreover, $dim(V_{3}) = dim(V_{4}) = 1$.
\end{theorem}
{\it Proof}. The variety $(V_{3})_{\sigma_{1}}$ is the variety $LZ$
of bands defined by $yx \approx y$. We obtain that: $(V_{3})_{\sigma_{1}}$
is different as $V_{3}$, therefore $(V_{3})_{\sigma_{1}}$ is proper
and $(V_{3})_{\sigma_{1}} \subseteq V_{3}$. Note, that the derived
variety $(V_{3})_{\sigma_{2}}$ is proper and is the variety $RZ$ of
right-zero semigroups but is not included in $V_{3}$. Similarly as
in the previous theorem  we conclude, that the variety $SL$ of
semilattices is not a derived variety of $V_{3}$, as all the
identities of $V_{3}$ are leftmost. In order to exclude that, the
variety $V_{1}$ defined by the identity (1) $zxy\approx zyx$ is  the
derived variety of $V_{3}$ by $\sigma_{5}$ consider the derived
identity $\sigma_{5}(zxy)\approx\sigma_{5}(zyx)$ of (1) by $\sigma_{5}$,
i.e. the identity $zxyxz\approx zyxyz$. If this identity would be satisfied
in $V_{3}$, then the identity $zxyz\approx zyxz$ would be satisfied in
$V_{3}$, which is not true due to the results of \cite{fen1}--\cite{3}.
Dually, the derived variety of $V_{3}$ by $\sigma_{6}$
is not the variety $V_{1}$. Therefore we conclude, that
$dim(V_{3})=1$. Similarly one can prove that $dim(V_{4})=1$. $\Box$

\begin{theorem}
The varieties $V_{5}$ and $V_{6}$ of bands defined by the identities:
\begin{center}
(5) $zxy \approx zxzy$ \hspace{5mm} and \hspace{5mm} (6) $yxz \approx yzxz$,
respectively,
\end{center}
are mutually derived by $\sigma_{4}$. Moreover, $dim(V_{5}) = dim(V_{6}) = 3$.
\end{theorem}
{\it Proof}. The proof that $(V_{5})_{\sigma_{1}}$
($(V_{5})_{\sigma_{2}}$ is the variety LZ (RZ) of left (right) zero
semigroups follows from the proof of previous observations.
Obviously: $\sigma_{4}(V_{5}) = V_{6}$, as the derived identity
$\sigma_{4}(yxz) \approx \sigma_{4}(yzxz)$ of (6) by $\sigma_{4}$
gives rise to the identity (5) $zxy\approx zxzy$ and vice versa.
Therefore $V_{6}$ and $V_{5}$ are mutually derived by $\sigma_{4}$.
We will show that the derived variety of $V_{5}$ by the
hypersubstitution $\sigma_{5}$ is the variety $V_{3}$, i.e.
$\sigma_{5}(V_{5}) = V_{3}$. To show this consider the derived
identity of (3) by $\sigma_{5}$, i.e. the identity $\sigma_{5}(yx)
\approx \sigma_{5}(yxy)$. This gives rise to the identity $yxy
\approx yxyxy$, which is obviously satisfied in $V_{5}$. Moreover,
note that the derived identity of (1) by $\sigma_{5}$, i.e.
$\sigma_{5}(zxy) \approx \sigma_{5}(zyx)$ gives rise to the identity
$zxyxz \approx zyxyz$, which can not be satisfied in $V_{5}$, as
otherwise the identity $zxyz \approx zyxz$ would be satisfied in
$V_{3}$, which is impossible by the results of \cite{fen1}--\cite{3}
and it has been shown already in the proof of theorem 4.5. Similarly
one can show, that the derived variety of $V_{5}$ by $\sigma_{6}$ is
the variety $V_{4}$. We conclude that $dim(V_{5}) = 3$. Similarly,
$dim(V_{6}) = 3$. $\Box$
\begin{definition}
The variety $NB$ of normal bands is defined by the identity:
\begin{center}
(7) $zxyz \approx zyxz$.
\end{center}
\end{definition}
\begin{theorem}
$dim(NB) = 4$.
\end{theorem}
{\it Proof}. For solidity of the variety $NB$ confront \cite[p. 96]{2}.
It follows, that all derived varieties of the variety of $NB$ are included
in the variety of $NB$. Similarly as before we show that $(NB)_{\sigma_{1}}$
is the variety $LZ$ of left-zero semigroups and $(NB)_{\sigma_{2}}$ is the
variety $RZ$ of right-zero semigroups. Both of them are proper subvarieties
of NB. It is obvious that $(NB)_{\sigma_{3}} = (NB)_{\sigma_{4}} = NB$.
We show only that $(NB)_{\sigma_{5}} = V_{1}$, as the derived indentity of
(1) $zxy \approx zyx$ by $\sigma_{5}$, i.e. $\sigma_{5}(zxy)\approx\sigma_{5}(zyx)$
gives rise to the identity $zxyxz\approx zyxyz$ satisfied in $NB$. In order to
exclude that the variety $LZ$ of left zero semigroups, defined by the
identity $yx\approx y$ equals to $(NB)_{\sigma_{5}}$, notice that the derived
identity of $yx \approx  y$ by $\sigma_{5}$ is the identity $yxy\approx y$, which is not satisfied
in $NB$, as the variety of $NB$ is defined by the set of regular identities
(cf. \cite{JP}), which has only regular consequences.
Similarly $(NB)_{\sigma_{6}} = V_{2}$, as the derived identity of
(2) $yxz \approx  xyz$ by $\sigma_{6}$, i.e.
$\sigma_{6}(yxz) \approx  \sigma_{6}(xyz)$ gives rise to the identity $zxyxz\approx zyxyz$
satisfied in $NB$ and we conclude that $dim(NB)=4$.  $\Box$

\begin{theorem}
$dim(RegB)=4$.
\end{theorem}
{\it Proof}. For solidity of the variety $RegB$ confront \cite[p.
96]{2}. Two derived subvarieties of $RegB$ are $LZ$ and $RZ$, by
$\sigma_{1}$ and $\sigma_{2}$, respectively. The derived varieties
of $RegB$ via $\sigma_{3}$ and $\sigma_{4}$ are equal to $RegB$. We
show that $(RegB)_{\sigma_{5}} = V_{3}$. To prove that, consider the
derived identity of the identity (3) $yx \approx yxy$ by
$\sigma_{5}$, i.e. the identity $\sigma_{5}(yx) \approx
\sigma_{5}(yxy)$ which gives rise to the identity $yxy \approx
yxyxy$ which is satisfied in $RegB$. In order to show that the
derived variety of $RegB$ by $\sigma_{5}$ is not equal to the
variety $V_{1}$, note that the derived identity of (1) $zxy \approx
zyx$ by $\sigma_{5}$, i.e. the identity $\sigma_{5}(zxy) \approx
\sigma_{5}(zyx)$ gives rise to the identity $zxyxz \approx zyxyz$
which is not satisfied in $RegB$, as it was shown in the proof of
theorem 4.6 that this identity is not satisfied in $V_{5}$, which is
a subvariety of $RegB$. Similarly, one can show, that the derived
variety of $RegB$ by $\sigma_{6}$ is the variety $V_{4}$. This
finishes the proof that $dim(RegB)=4$. $\Box$

We expressed the situation of theorems above on the diagram:

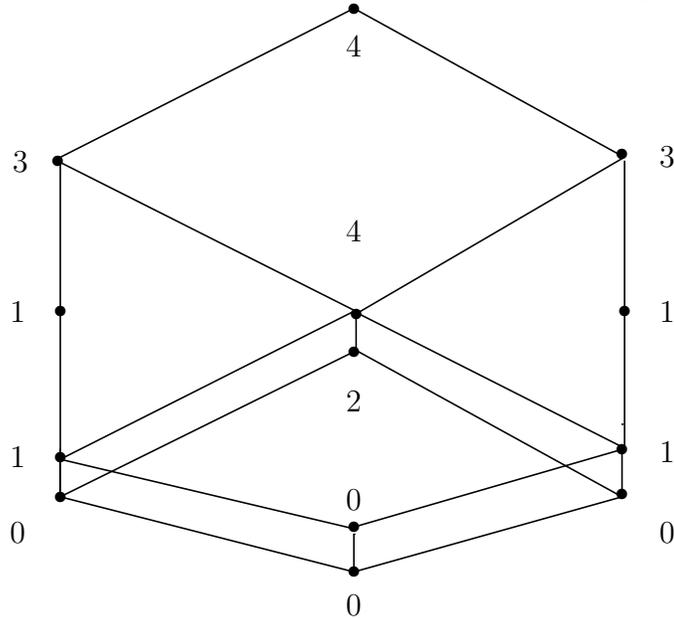
\vspace{15mm}

\input{dim2.pic}

{\bf Acknowledgement}. The first author expresses her thanks to the referee
for his valuable comments.

\end{document}